\documentclass[12pt,a4paper,reqno]{amsart}

\usepackage{amssymb}
\usepackage{amscd}
\usepackage{amsfonts}
\usepackage{setspace}

\numberwithin{equation}{section}
\theoremstyle{plain}
\newtheorem{theorem}[subsection]{Theorem}
\newtheorem{proposition}[subsection]{Proposition}
\newtheorem{lemma}[subsection]{Lemma}
\newtheorem{corollary}[subsection]{Corollary}

\theoremstyle{definition}

\newtheorem{definition}[subsection]{Definition}

\newtheorem*{mainthm1-again}{Theorem \ref{mainthm1}'}
\newtheorem*{mainthm2-again}{Theorem \ref{mainthm2}'}
\newtheorem*{frei-ruz-again}{Theorem \ref{frei-ruz}'}
\newtheorem*{key-prop-2-repeat}{Proposition \ref{key-prop-2}}
\renewcommand{\leq}{\leqslant}
\renewcommand{\geq}{\geqslant}
\newsavebox{\proofbox}
\savebox{\proofbox}{\begin{picture}(7,7)  \put(0,0){\framebox(7,7){}}\end{picture}}

\newcommand\id{\operatorname{id}}

\newcommand\alg{\operatorname{alg}}
\newcommand\gp{\operatorname{gp}}

\onehalfspace
\begin{document}
\title{Approximate groups, I:  the torsion-free nilpotent case}
\author{Emmanuel Breuillard}
\address{Laboratoire de Math\'ematiques 
Universit\'e Paris-Sud 11, 91405 Orsay cedex,
France }
\email{emmanuel.breuillard@math.u-psud.fr}
\author{Ben Green}
\address{Centre for Mathematical Sciences\\
Wilberforce Road\\
Cambridge CB3 0WA\\
England }
\email{b.j.green@dpmms.cam.ac.uk}
\subjclass{}

\begin{abstract}
We describe the structure of ``$K$-approximate subgroups'' of torsion-free
nilpotent groups, paying particular attention to Lie groups.

Three other works, by Fisher-Katz-Peng, Sanders and Tao, have appeared which
independently address related issues. We comment briefly on some of the
connections between these papers.
\end{abstract}

\maketitle
\tableofcontents

\setcounter{tocdepth}{1}

\section{Introduction}

\textsc{Approximate groups.} A fair proportion of the subject of additive
combinatorics is concerned with approximate analogues of exact algebraic
properties, and the extent to which they resemble those algebraic
properties. In this paper we are concerned with approximate groups.

By an \emph{ambient group} we simply mean some group in which all the
objects being discussed are contained, so that it makes sense to talk about
multiplication of elements, inverses and the identity element. Suppose that $
A$ is a finite set in some ambient group. What does it mean to say that $A$
is an approximate subgroup?

It is well-known to all students of group theory that $A$ is a genuine
subgroup if, and only if, we have $xy^{-1} \in A$ whenever $x,y \in A$.
Perhaps the most natural way in which a set $A$ may be \emph{approximately}
a subgroup, then, is if the set $A A^{-1} = \{xy^{-1} : x,y \in A\}$ has
cardinality not much bigger than $|A|$, perhaps $|AA^{-1}| \leq K|A|$ for
some constant $K$.

Sets with this property are said to have \emph{small doubling} and this is
indeed a commonly encountered condition in additive combinatorics. It is a
perfectly workable notion of approximate group in the abelian setting and
the celebrated Fre\u{\i}man-Ruzsa theorem describes subsets of $\mathbb{Z}$
with this property (we will state it below). However in the foundational
work of Tao \cite{tao-noncommutative} it was noted that in noncommutative
settings a somewhat different, though closely related, notion of approximate
group is more natural. We now give Tao's definition.

\begin{definition}[Approximate groups]
Let $K \geq 1$. A set $A$ in some ambient group is called a $K$-approximate
group if

\begin{enumerate}
\item It is symmetric, i.e. if $a \in A$ then $a^{-1} \in A$, and the
identity lies in $A$;

\item There is a symmetric subset $X$ lying in $A\cdot A$ with $|X| \leq K$
such that $A \cdot A \subseteq X \cdot A$.
\end{enumerate}
\end{definition}

This definition gives rise to some very pleasant properties, and we shall
list them in \S \ref{basic-add-sec}. In that section we also briefly recall
the relation between approximate groups in this sense and sets with small
doubling.

Our aim in this paper is to ``describe'' the structure of approximate
subgroups of torsion-free nilpotent groups in terms of more explicit
algebraic objects. A companion paper \cite{breuillard-green-solv} tackles
the same question for solvable subgroups of $\mbox{GL}_d(\mathbb{C})$. Tao 
\cite{tao-solvable} has addressed questions of this type, working in fact
with solvable groups in general. In his paper he introduces the following
rather nice paradigm for ``describing'' sets by others.

\begin{definition}[Control]
Suppose that $A$ and $B$ are two sets in some ambient group, and that $K
\geq 1$ is a parameter. We say that $A$ is $K$-controlled by $B$, or that $B$
$K$-controls $A$, if $|B| \leq K|A|$ and there is some set $X$ in the
ambient group with $|X| \leq K$ and such that $A \subseteq (X \cdot B) \cap
(B \cdot X)$.
\end{definition}

This is essentially equivalent to saying that $A$ and $B$ have roughly the
same size and that $A$ is covered by a few left-translates of $B$, and also
by a few right translates of $B$. Indeed if $A \subseteq \bigcup_{i=1}^k x_i
B$ and also $A \subseteq \bigcup_{j=1}^{l} B y_j$ then we may take $X =
\{x_1,\dots,x_k,y_1,\dots,y_l\}$ in the definition above; the other
direction of the equivalence is even more obvious.

In \S \ref{basic-add-sec} we will discuss (following Tao's paper extremely
closely) how this notion of control interacts with the aforementioned
notions of approximate group and small doubling.\vspace{11pt}

The structure of approximate subgroups of torsion-free \emph{abelian} groups
is described by the Fre\u{\i}man-Ruzsa theorem \cite{freiman-book,ruzsa-freiman}. The bounds in the following version of it,
which is stated in the language introduced above, are due to Chang \cite{chang-freiman}. Here and for the remainder of the paper the letter $C$
represents an absolute constant which could be computed explicitly if
desired: different instances of the letter may denote different constants.
We will often use subscripts to indicate dependence on other parameters: for
example, $C_{s}$ is an absolute constant depending on $s$.

\begin{theorem}[Fre{\u{\i}}man-Ruzsa; Chang]
\label{frei-ruz} Let $G$ be a torsion-free abelian group and let $K \geq 1$
be a parameter. Suppose that $A \subseteq \mathbb{Z}$ is a $K$-approximate
group. Then $A$ is $e^{CK^C}$-controlled by a set $P$ of the form 
\begin{equation*}
P = \{l_1 x_1 + \dots + l_k x_k : |l_1| \leq L_1, \dots , |l_k| \leq L_k\},
\end{equation*}
for some $x_1\dots,x_k \in \mathbb{Z}$, where $k \leq CK^C$.
\end{theorem}

A set $P$ of this form is called a \emph{generalised arithmetic progression}
, or progression for short\footnote{There is a very slight difference between the terminology used in this paper
and that which is standard: typically, a progression is a set of the form $\{x_0 + l_1 x_1 + \dots + l_k x_k : 0 \leq l_i < L_i\}$. We have found it
convenient to disallow the presence of $x_0$ in this paper and to use the
more symmetric condition $|l_i| \leq L_i$. Note however that every
progression in our sense is economically contained in one according to the
more standard definition, and vice versa (though one might need to increase
the dimension by 1).}. The number $k$ is referred to as the \emph{dimension}
of $k$.

The Fre\u{\i}man-Ruzsa theorem is usually stated and proved only for subsets
of $\mathbb{Z}$ and not for torsion-free abelian groups in general. Simple
modifications allow one to obtain the more general statement, and we will
remark further on this later on.

Let us turn now to nilpotent groups, pausing to recall the definition. Let $G
$ be a group and suppose that $s \geq 1$ is an integer. If the lower central
series defined by 
\begin{equation*}
G_0 = G_1 = G, G_2 = [G,G_1], G_3 = [G,G_2],\dots
\end{equation*}
terminates with $G_{s+1} = \{\id_G\}$ then we say that $G$ is $s$-step
nilpotent. A prototypical example of a torsion-free nilpotent group $G$ is a
group of upper triangular matrices with ones on the diagonal, such as the
Heisenberg group $G = \left( 
\begin{smallmatrix}
1 & \mathbb{R} & \mathbb{R} \\ 
0 & 1 & \mathbb{R} \\ 
0 & 0 & 1
\end{smallmatrix}
\right)$, which is an example of a $2$-step nilpotent (Lie) group.

Here is an example of an approximate subgroup of the Heisenberg group. It is
also discussed quite explicitly in Tao's paper \cite{tao-solvable}.

Take 
\begin{equation*}
u_1 = \left( 
\begin{smallmatrix}
1 & 1 & 0 \\ 
0 & 1 & 0 \\ 
0 & 0 & 1
\end{smallmatrix}
\right) \quad \mbox{and} \quad u_2 = \left( 
\begin{smallmatrix}
1 & 0 & 0 \\ 
0 & 1 & 1 \\ 
0 & 0 & 1
\end{smallmatrix}
\right),
\end{equation*}
and consider also the commutator 
\begin{equation*}
[u_1,u_2] := u_1^{-1}u_2^{-1}u_1 u_2 = \left( 
\begin{smallmatrix}
1 & 0 & 1 \\ 
0 & 1 & 0 \\ 
0 & 0 & 1
\end{smallmatrix}
\right).
\end{equation*}
Let $L_1,L_2 \geq 1$ be integers. Then the set $A :=
\{u_1^{l_1}u_2^{l_2}[u_1,u_2]^{l_{12}} : |l_1| \leq L_1, |l_2| \leq L_2,
|l_{12}| \leq L_{1}L_{2}\}$ consists of the matrices 
\begin{equation*}
\{\left( 
\begin{smallmatrix}
1 & x & z \\ 
0 & 1 & y \\ 
0 & 0 & 1
\end{smallmatrix}
\right) : |x| \leq L_1, |y| \leq L_2, |z| \leq L_1L_2\}
\end{equation*}
and it is not hard to check that $A \cup A^{-1}$ is a $K$-approximate group
for some absolute constant $K$.

The set $A$ is obviously a close analogue of the generalised progressions
considered in the abelian setting. The construction may be generalised, but
to do this we must first discuss commutators in more detail. Suppose that $G$
is an $s$-step nilpotent group and that $u_{1},\dots ,u_{k}\in G$. We may
inductively assign a \emph{weight vector} $\chi \in \mathbb{N}_{0}^{k}$ to
every (formal) commutator involving the $u_{i}$s by setting $\chi
(u_{i})=e_{i}$ and defining $\chi $ inductively on higher commutators via $
\chi ([c,c^{\prime }])=\chi (c)+\chi (c^{\prime })$. Thus if $k=4$ then $
\chi ([u_{1},u_{2}])=(1,1,0,0)$ and $\chi ([u_{1},[u_{2},u_{4}])=(1,1,0,1)$.
We follow M.~Hall \cite[Chapter 11]{Hall} in defining \emph{basic commutators
}. This is a (non-unique) extension of $u_{1},\dots ,u_{k}$ to an ordered
list $u_{1}\prec \dots \prec u_{t}$ in which $u_{k+1},\dots $ are certain
commutators involving $u_{1},\dots ,u_{k}$. We suppose that they are ordered
so that commutators with the same weight vector are consecutive, and so that
higher order commutators come before lower order ones. If $c_{i},c_{j}$ have
already been admitted as basic commutators then $c_{k}=[c_{i},c_{j}]$
qualifies as basic if $c_{i}\succ c_{j}$ and if, writing $c_{i}=[c_{s},c_{t}]
$, $c_{j}\succeq c_{t}$.

For example when $k=s=3$ we have $t=14$, a possible listing of the basic
commutators being $u_{1}$, $u_{2}$, $u_{3}$, $[u_{2},u_{1}]$, $[u_{3},u_{2}]$
, $[u_{3},u_{1}]$, $[[u_{2},u_{1}],u_{1}]$, $[[u_{2},u_{1}],u_{2}]$, $
[[u_{2},u_{1}],u_{3}]$, $[[u_{3},u_{1}],u_{1}]$, $[[u_{3},u_{1}],u_{2}]$, $
[[u_{3},u_{1}],u_{3}]$, $[[u_{3},u_{2}],u_{2}]$, $[[u_{3},u_{2}],u_{3}]$.
Note incidentally the formula of Witt, which states that the number of basic
commutators of order $r$ on $k$ generators is $\frac{1}{r}\sum_{d|r}\mu
(d)k^{r/d}$.

Write $\chi(j)$ for the weight vector of the commutator $u_j$. If $L =
(L_1,\dots,L_k)$ is a vector of positive integers and $\chi \in \mathbb{N}
_0^k$, we define $L^{\chi} := L_1^{\chi_1}\dots L_k^{\chi_k}$.

\begin{definition}[Nilpotent progressions]
Suppose that $G$ is an $s$-step nilpotent group and that $u_1,\dots,u_k \in G
$. Let $L = (L_1,\dots,L_k)$ be a vector of positive integers. Then the
nilpotent progression $P(u_1,\dots,u_k;L)$ on generators $u_1,\dots,u_k$
with lengths $L$ is the set $\{u_1^{l_1} \dots u_t^{l_t} : |l_j| \leq
L^{\chi(j)}\}$, where $u_1,u_2,\dots,u_t$ is the ordered list of basic
commutators involving the $u_i$.
\end{definition}

Any $s$-step nilpotent progression is the homomorphic image of a nilpotent
progression in $\Gamma_{k,s}$, the \emph{free} $s$-step nilpotent group on $k
$ generators. We have found this to be the right way to think about
nilpotent progressions: to study them, one should establish
homomorphism-invariant properties of nilpotent progressions in the free
nilpotent group.

Consideration of the free case reveals our reason for involving only \emph{basic} commutators, since by restricting to these the elements $
u_1^{l_1}\dots u_t^{l_t}$ of a nilpotent progression in the free nilpotent
group are all distinct. This follows from the results of \cite[Chapter 11]
{Hall}.\vspace{11pt}

We are now in a position to state our main theorem, which is the analogue of
the Fre\u{\i}man-Ruzsa theorem in the nilpotent setting.

\begin{theorem}
\label{mainthm0} Let $\Gamma $ be a torsion-free $s$-step nilpotent group,
and suppose that $A\subseteq \Gamma $ is a $K$-approximate subgroup. Then
there are elements $u_{1},\dots ,u_{k}\in \Gamma $, $k\leq K^{C_{s}}$, and
lengths $L=(L_{1},\dots ,L_{k})$ such that $A$ is $e^{K^{C_{s}}}$-controlled
by the nilpotent progression $P(u_{1},\dots ,u_{k};L)$.
\end{theorem}

The proof of this theorem occupies the majority of the paper. However in the
later sections we gather some properties of nilpotent progressions which may
be of interest in their own right. We also comment on the connection between
our results and those of Sanders \cite{sanders}.\vspace{11pt}

\textsc{acknowledgement.} It is a pleasure to thank Elon Lindenstrauss, Tom Sanders and Terry Tao for a number of helpful conversations.

\section{Strategy of the proof}

\label{strategy-sec}

The key idea for establishing Theorem \ref{mainthm0}, already implicit in 
\cite{tao-noncommutative} and described explicitly in \cite{fkp}, is to use
a little Lie theory. We may clearly suppose, in proving Theorem \ref
{mainthm0}, that $\Gamma $ is finitely-generated. An embedding theorem of
Mal'cev \cite{Rag} states that every finitely-generated torsion-free
nilpotent group embeds as a co-compact discrete subgroup of a
simply-connected nilpotent Lie group of the same step. It therefore suffices
to establish Theorem \ref{mainthm0} when $\Gamma $ is a subgroup of a
simply-connected $s$-step nilpotent Lie group, say $G$ (the Heisenberg group
is an example of such a group).

Working in this setting enables us to exploit the Lie algebra $\mathfrak{g}
=\log G.$ It is well-known in the theory of simply-connected nilpotent Lie
groups (see, for example, \cite{bourbaki}) that there are mutually inverse
diffeomorphisms $\exp :\mathfrak{g}\rightarrow G$ and $\log :G\rightarrow 
\mathfrak{g}$ between the group $G$ and its Lie algebra $\mathfrak{g}$,
which is a vector space $\mathbb{R}^{\dim (G)}$ together with an additional
bracket operation $[,]:\mathfrak{g}\times \mathfrak{g}\rightarrow \mathfrak{g
}$ which is antisymmetric, bilinear and satisfies the Jacobi identity.

In the case of the Heisenberg group we may identify $\mathfrak{g}$ with the
vector space 
\begin{equation*}
\left( 
\begin{smallmatrix}
0 & x & z \\ 
0 & 0 & y \\ 
0 & 0 & 0
\end{smallmatrix}
\right).
\end{equation*}
The exponential map $\exp : \mathfrak{g} \rightarrow G$ is then simply the
usual exponentiation of matrices.

In the Lie algebra setting it is quite natural to consider a different type
of nilpotent progression which, to distinguish it from the nilpotent
progressions already described, we call a \emph{nilbox}. To define nilboxes
we must first describe commutators in $\mathfrak{g}$.

We consider first the free nilpotent Lie algebra $\mathfrak{n}_{k,s}$ with
generators $X_1,\dots,X_k$. We will be looking at higher order commutators
such as $[X_1,[X_2,X_5]]$, and once again we will associate a weight vector $
\chi \in \mathbb{N}_0^k$ to each of these. The definition is the same as
before (in the example just given, $\chi = (1,1,0,0,1,\dots))$. We have a
decomposition $\mathfrak{n}_{k,s} = \oplus_{\chi} V_{\chi}$ into \emph{
weight spaces} $V_{\chi}$, where $V_\chi$ consists of commutators with a
fixed weight $\chi$. Just as for group commutators, we may extend $
X_1,\dots,X_k$ to an ordered list $X_1,\dots,X_t$ of \emph{basic commutators}
. The definition of these is precisely the same as for group commutators,
except that the bracket now refers to the Lie algebra operation rather than
the group commutator.

The reason for introducing \emph{basic} commutators becomes clear in this
context: by a theorem of Witt \cite[Chapter 11]{Hall} the elements $
X_1,\dots,X_t$ form a basis for $\mathfrak{n}_{k,s}$ as a vector space over $
\mathbb{C}$. We call this an \emph{adapted basis} for $\mathfrak{n}_{k,s}$.

\begin{definition}[Free nilboxes]
Let $k,s \geq 1$ be integers, and let $X_1,\dots,X_k$ be generators for the
free nilpotent Lie algebra $\mathfrak{n}_{k,s}$. Let $L = (L_1,\dots,L_k)$
be a vector of positive integer lengths. Then the free $s$-step nilbox with
lengths $L$ is the set $\mathfrak{B}(X_1,\dots,X_k;L) \subseteq \mathfrak{n}
_{k,s}$ defined by 
\begin{equation*}
\mathfrak{B}(X_1,\dots,X_k;L) = \{l_1 X_1 + \dots + l_t X_t : |l_j| \leq
L^{\chi(j)}\},
\end{equation*}
where $X_1,\dots,X_t$ is an adapted basis for $\mathfrak{n}_{k,s}$.
\end{definition}

\begin{definition}[Nilboxes]
Let $\mathfrak{g}$ be an $s$-step nilpotent Lie algebra and suppose that $
x_1,\dots,x_k \in \mathfrak{g}$. Let $L = (L_1,\dots,L_k)$ be a vector of
positive integer lengths. Then we define the nilbox $\mathfrak{B}
(x_1,\dots,x_k;L)$ to be the image $\pi(\mathfrak{B}(X_1,\dots,X_k;L))$,
where $\pi : \mathfrak{n}_{k,s} \rightarrow \mathfrak{g}$ is the Lie algebra
homomorphism induced by mapping $X_i$ to $x_i$, $i=1,\dots,k$. It is
convenient to write $x_i := \pi(X_i)$ for $i = k+1,\dots,t$ also.
\end{definition}

As we remarked, it suffices to prove Theorem \ref{mainthm0} when $\Gamma $
is a subgroup of a simply-connected $s$-step nilpotent Lie group $G$. We may
now divide this task into the task of proving the following two propositions.

\begin{proposition}[Control by nilboxes]
\label{key-prop-1} Suppose that $G$ is an $s$-step simply-connected
nilpotent Lie group with Lie algebra $\mathfrak{g}$, and that $A \subseteq G$
is a $K$-approximate group. Then there are $x_1,\dots,x_k \in \mathfrak{g}$
such that

\begin{enumerate}
\item $k \leq K^{C_s}$;

\item $\exp(x_1),\dots, \exp(x_t)$ lie in the group $\langle A\rangle$
generated by $A$;

\item There is a nilbox $\mathfrak{B}(x_1,\dots,x_k;L)$ such that $\exp(
\mathfrak{B}(x_1,\dots,x_k;L))$ $e^{K^{C_s}}$-controls $A$.
\end{enumerate}
\end{proposition}

\begin{proposition}[Nilpotent progressions control nilboxes]
\label{key-prop-2} Suppose that $G$ is an $s$-step simply-connected
nilpotent Lie group with Lie algebra $\mathfrak{g}$. Suppose that $
x_{1},\dots ,x_{k}\in G$, and write $u_{i}:=\exp (x_{i})$. Let $
L=(L_{1},\dots ,L_{k})$ be a vector of positive integer lengths. Then the
nilpotent progression $P(u_{1},\dots ,u_{k};L)$ $e^{k^{C_{s}}}$-controls $
\exp (\mathfrak{B}(x_{1},\dots ,x_{k};L))$.
\end{proposition}

Theorem \ref{mainthm0} clearly follows from the combination of the last two
propositions after we observe that the nilpotent progression $P(u_{1},\dots
,u_{k};L)$ obtained in Proposition \ref{key-prop-2} entirely lies in $\Gamma 
$ and that if $A$ and $B$ are two subsets of the subgroup $\Gamma $ and 
$A$ is $K$-controlled by $B$ in $G,$ then $A$ must also 
be $K$-controlled by $B$ in $\Gamma$.

Either proposition is conceivably of independent interest. For example it
seems to be easier to study nilboxes than nilpotent progressions. The proof
of Proposition \ref{key-prop-1} is essentially additive-combinatorial and
occupies the next four sections. The proof of Proposition \ref{key-prop-2}
requires a certain amount of material on coordinates in nilpotent Lie
groups: this material is summarised in Appendix \ref{coord-append} and the
proposition itself is confirmed in Section \ref{nilbox}.

To conclude this section let us note that there is a certain arbitrariness
in the definitions of nilpotent progression and nilbox, coming from the
noncanonical choice of an ordering for the basic commutators (or indeed for
the basic commutators themselves, which are defined in different ways by
different authors). This is not a serious matter and any other choice would
lead to completely equivalent theorems.

\section{Preliminaries from multiplicative combinatorics}

\label{basic-add-sec}

We take the opportunity to record some basic facts about noncommutative
product sets, and in particular concerning the notions of $K$-approximate
group and $K$-control defined in the introduction. This material is all due
to Tao \cite{tao-noncommutative}; in turn some of that is based on earlier
work of Ruzsa in the abelian setting. See also the book of Tao and Vu \cite
{tao-vu}, especially \S 2.7.

\begin{proposition}[Approximate groups and control]
\label{approx-gp-basics} Let $K \geq 1$ be a parameter and let $A$ be a set
in some ambient group $G$. If $n \geq 1$ is an integer we write $A^n = \{a_1
\dots a_n : a_1,\dots,a_n \in A\}$ and $A^{\pm n} = \{a_1^{\varepsilon_1}
\dots a_n^{\varepsilon_n} : a_1,\dots,a_n \in A,
\varepsilon_1,\dots,\varepsilon_n \in \{-1,1\}\}$.

\begin{enumerate}
\item If $\pi : G \rightarrow H$ is a homomorphism and if $A \subseteq G$ is
a $K$-approximate group then $\pi(A)$ is a $K$-approximate subgroup of $H$.

\item If $A$ is a $K$-approximate group then $|A^{\pm n}| = |A^n| \leq
K^{n-1}|A|$ and $A^n$ is $K^{n+1}$-controlled by $A$.

\item If $B,C$ are further subsets of $G$ and if $A$ is $K$-controlled by $B$
and $B$ is $K$-controlled by $C$, then $A$ is $K^{2}$-controlled by $C$.

\item If $A$ and $B$ are $K$-approximate groups and $A$ is $K$-controlled by 
$B,$ then $B$ is $K^{4}$-controlled by $A.$

\item If the doubling constant $|A^{2}|/|A|$ is at most $K$ then there is an 
$f_{1}(K)$-approximate group $B\subseteq A^{\pm 3}$ which $f_{2}(K)$
-controls $A$. If the tripling constant $|A^{3}|/|A|$ is at most $K$ then we
may take $B=A^{\pm 3}$.

\item If $A$ is a $K$-approximate group and if $A^{\prime }\subseteq A$ is a
subset with $|A^{\prime }|\geq |A|/K$ then $A^{\prime \pm 3}$ is an $
f_{3}(K) $-approximate group which $f_{4}(K)$-controls $A$. The same is in
fact true under the essentially weaker assumption that $|A^{3}|\leq K|A|$.
\end{enumerate}

All of the quantities $f_1(K),\dots,f_4(K)$ can be taken to be polynomial in 
$K$.
\end{proposition}

\noindent\textit{Proof. } (i) follows immediately from the definition. To
prove (ii), suppose that $X$ is a symmetric set such that $|X| \leq K$ and $
A\cdot A \subseteq X\cdot A$. Then, since $A$ is symmetric, we also have $
A\cdot A \subseteq A \cdot X$ and hence $A^n \subseteq X^{n-1}\cdot A$ and $
A^n \subseteq A\cdot X^{n-1}$ for all $n \geq 1$ by an easy induction, from
which the result follows immediately.

Parts (iii) is very easy and follows straight from the definitions. Part
(iv) follows from the non-abelian Ruzsa covering lemma (\cite
{tao-noncommutative} Lemma 3.6). Parts (v) and (vi) may be found in \cite
{tao-noncommutative}: (v) is Theorems 4.6 and Corollary 3.10 of that paper,
whilst (vi) follows from Lemma 3.6 and Corollary 3.10.\hfill {\ 
\usebox{\proofbox}}\vspace{11pt}

\section{Generalized arithmetic progressions in Lie algebras}

\label{prog-section}

Let $\mathfrak{g}$ be an $s$-step nilpotent Lie algebra and let $\mathfrak{p}
\subseteq \mathfrak{g}$ be a progression, thus 
\begin{equation*}
\mathfrak{p} = \{l_1 x_1 + \dots + l_k x_k : |l_i| \leq L_i\}
\end{equation*}
for some lengths $L_1,\dots,L_k$ and some $x_1,\dots,x_k \in \mathfrak{g}$.
Our aim is to understand ways in which $\mathfrak{p}$ can interact with the
bracket operation $[,]$.

We begin with a definition.

\begin{definition}[Nilcompletion]
\label{nilcompletion-def} Suppose that $\mathfrak{b}\subseteq \mathfrak{g}$
is a set. Then by the \emph{nilcompletion} $\overline{\mathfrak{b}}$ of $
\mathfrak{b}$ we mean the set $\mathfrak{b}+[\mathfrak{b},\mathfrak{b}]+[ 
\mathfrak{b},[\mathfrak{b},\mathfrak{b}]]+\dots +[[\mathfrak{b},\mathfrak{b}
],[\mathfrak{b},[\mathfrak{b},\mathfrak{b}]]]\dots $, where the sum is over
all\footnote{
In this paper we are not concerned with the dependence of our estimates on
the step parameter $s$. If we were, it might be more efficient at this point
to take only \emph{nested} commutators into account in making this
definition.} commutators.
\end{definition}

Now it is known from standard Lie theory that if $\mathfrak{g}_{i}:=\log
(G_{i})$ then $[\mathfrak{g}_{i},\mathfrak{g}_{j}]\subseteq \mathfrak{g}
_{i+j}$, and so any commutator with more than $s$ copies of $\mathfrak{b}$
vanishes identically. The number of commutators of order $k+1$ is the $k$-th
Catalan number $C_{k}=\frac{1}{k+1}\binom{2k}{k}$. We easily see that the
total number of such up to order $s$ can be bounded above by $4^{s}.$

\begin{lemma}[Properties of the nilcompletion]
\label{lem5.1} Let $\mathfrak{b} \subseteq \mathfrak{g}$ be a set. Then:

\begin{enumerate}
\item For any integer $m \geq 1$ we have $\overline{m\mathfrak{b}} \subseteq
m^s\overline{\mathfrak{b}}$;

\item $[\overline{\mathfrak{b}},\overline{\mathfrak{b}}] \subseteq \overline{
\mathfrak{b}}$.
\end{enumerate}
\end{lemma}

\noindent \textit{Proof. } The first inclusion is a consequence of the fact
that an $s$-fold commutator is $s$-multilinear. For example (when $m=2$) $
[b_{1}+b_{1}^{\prime },b_{2}+b_{2}^{\prime }]$ may be written as a sum of
four commutators involving only elements of $\mathfrak{b}$. The second is
also immediate, and it is the main reason for introducing this definition of
nilcompletion.\hfill {\usebox{\proofbox}}\vspace{11pt}

Now if $\mathfrak{p}$ is a progression then its nilcompletion $\overline{
\mathfrak{p} }$ need not be in general. Proposition \ref{prop33} below is a
good substitute for the failure of this statement, however. Before stating
it we record a simple number-theoretic lemma that we will use several times
in the sequel.

\begin{lemma}
\label{num-th-lemma} Suppose that $L_1,\dots,L_k$ are integers. Then every
nonnegative integer less than or equal to $L_1\dots L_k$ can be written as
the sum of at most $2^{k-1}$ numbers of the form $l_1\dots l_k$ with $1 \leq
l_i \leq L_i$ for $i=1,\dots,k$.
\end{lemma}

\noindent\textit{Proof. } This may be established by induction on $k$ from
the base case $k=2$: to prove that case, write a given $m$ as $qL_{2}+r$,
where $0\leq q<L_{1}$ and $1\leq r\leq L_{2}$.\hfill{\usebox{\proofbox}}
\vspace{11pt}

\begin{proposition}
\label{prop33} Suppose that $\mathfrak{p}\subseteq \mathfrak{g}$ is a 
\textup{progression} of dimension $k$ and that $\overline{\mathfrak{p}}$ is
its nilcompletion. Then there is a \textup{progression} $\mathfrak{q}$ of
dimension at most $(4k)^{s}$ such that $[\mathfrak{q},\mathfrak{q}]\subseteq
(2k)^{2s} \mathfrak{q}$ and $\overline{\mathfrak{p}}\subseteq \mathfrak{q}
\subseteq (2k)^{s}\overline{\mathfrak{p}}$.
\end{proposition}

\noindent \textit{Proof. } Suppose as usual that $\mathfrak{p}$ has
generators $x_{1},\dots ,x_{k}$. Then 
\begin{equation}
\mathfrak{b}\subseteq \mathfrak{p}\subseteq k\mathfrak{b},  \label{inc-1}
\end{equation}
where $\mathfrak{b}:=\bigcup_{i=1}^{k}\{l_{i}x_{i}:|l_{i}|\leq L_{i}\}$. Let 
$\mathcal{C}$ be the set of all commutators in the $x_{i}$'s, such as $
[[x_{1},x_{7}],[x_{3},x_{5}]]$. To each $c\in \mathcal{C}$ we may assign a
weight vector $\chi (c)\in \mathbb{N}_{0}^{k}$ much as described in the
introduction; for this example, $\chi =(1,0,1,0,1,0,1,0,0,\dots )$. We first
claim that 
\begin{equation}
\overline{\mathfrak{b}}\subseteq \mathfrak{r}\subseteq 2^{s-1}\overline{ 
\mathfrak{b}}.  \label{inc-3}
\end{equation}
where $\mathfrak{r}$ is the progression defined by $\mathfrak{r}
=\{\sum_{c\in \mathcal{C}}l_{c}c:|l_{c}|\leq L^{\chi (c)}\}$, where $L^{\chi
(s)}:=\prod_{i=1}^{k}L_{i}^{(\chi (s))_{i}}.$

The first inclusion is clear. To see the second, we invoke Lemma \ref
{num-th-lemma}. Now $\overline{ \mathfrak{b}}$ contains all elements of the
form $l^{\chi (c)}c$, provided $|l_{i}|\leq L_{i}$ for all $i=1,...,k$.
Hence by our observation we must have $\mathfrak{r}\subseteq n\overline{
\mathfrak{b}}$ where $n=2^{s-1}$. This proves (\ref{inc-3}).

Now from \eqref{inc-1}, \eqref{inc-3} and Lemma \ref{lem5.1} (i) we have $
\overline{\mathfrak{p}}\subseteq k^{s}\overline{\mathfrak{b}}\subseteq k^{s}
\mathfrak{r}$. This last set $k^{s}\mathfrak{r}$ is another progression, and
this we take to be $\mathfrak{q}$. Specifically, $\mathfrak{q}=k^{s}
\mathfrak{r}=\{\sum_{c\in \mathcal{C}}l_{c}c:|l_{c}|\leq k^{s}L^{\chi (c)}\}$
. The dimension of $\mathfrak{q}$ is $|\mathcal{C}|$, which is $\leq (4k)^{s}
$, since as we observed earlier, the total number of bracketing patterns of
length at most $s$ is $\leq 4^{s}$ . Now by construction we have $\overline{
\mathfrak{p}}\subseteq \mathfrak{q}$. Finally \eqref{inc-1} and \eqref{inc-3}
yield $\mathfrak{q}=k^{s}\mathfrak{r}\subseteq (2k)^{s}\overline{\mathfrak{b}
}\subseteq (2k)^{s}\overline{\mathfrak{p}}$, as required.

It remains to show that $[\mathfrak{q},\mathfrak{q}]$ is contained in an
appropriate multiple of $\mathfrak{q}$. This now follows from the fact that $
[m\mathfrak{b},m\mathfrak{b}]\subseteq m^{2}[\mathfrak{b},\mathfrak{b}]$ for
any set $\mathfrak{b}\subseteq \mathfrak{g}$ and from Lemma \ref{lem5.1}
.\hfill{\usebox{\proofbox}}\vspace{11pt}

\section{Some nilpotent algebra}

\label{sec4}

In this section we establish some results about nilpotent Lie groups and
their interaction with their Lie algebras via the exponential map.
Throughout this section, then, $G$ is a simply-connected $s$-step nilpotent
Lie group.

One of the key tools in this paper is a theorem of Lazard \cite{Laz1,Laz2},
given in Lemma \ref{stewart-lem} below, stating that both addition and the
bracket operation on $\mathfrak{g}$ may be expressed using words in $G$ of
length $C_{s}$. This consequence of the Baker-Campbell-Hausdorff formula was
discussed in detail in the 1969 thesis of Ian Stewart (see also \cite
{stewart}). Results of this type were first exploited in the
additive-combinatorial setting by Fisher, Katz and Peng \cite{fkp}.

We shall also need the following lemma on rational powers in nilpotent Lie
groups.

\begin{lemma}[Rational powers of products]
\label{rat-powers} Suppose that $G$ is an $s$-step nilpotent Lie group and
that $x_{1},\dots ,x_{t}\in G$. Let $\eta \in \mathbb{Q}$. Then there is an
integer $t^{\prime }$ and for $j\in \{1,...,t^{\prime }\}$ there are
polynomials $P_{j}(\eta )\in \mathbb{Q}[\eta ]$ drawn from a finite list
depending only on $s$ such that $(x_{1}\dots x_{t})^{\eta
}=x_{i_{1}}^{P_{1}(\eta )}\dots x_{i_{t^{\prime }}}^{P_{t^{\prime }}(\eta )}$
for some indices $i_{j}\in \{1,...,t\}.$
\end{lemma}

\noindent \textit{Proof. } This follows from the Hall-Petresco formula as
presented, for example, in \cite[Chapter 6]{serre} or \cite[Appendix A]
{du-sautoy}. This formula states (for an arbitrary $s$-step nilpotent group)
that there are words $w_{2}(x_{1},x_{2},\dots ,x_{t}),\dots
,w_{s}(x_{1},x_{2},\dots ,x_{t})$ such that $w_{j}\in G_{j}$ and 
\begin{equation*}
(x_{1}\dots x_{t})^{n}=x_{1}^{n}x_{2}^{n}\dots x_{t}^{n}w_{2}(x_{1},\dots
,x_{t})^{\binom{n}{2}}\dots w_{s}(x_{1},\dots ,x_{t})^{\binom{n}{s}}
\end{equation*}
for all positive integers $n$. By induction on the step of $G,$ this implies
that there are polynomials $P_{j}(n)$ with coefficients in $\mathbb{Q}$
depending only on $s$ such that $(x_{1}\dots x_{t})^{n}=\prod_{1\leq j\leq
t^{\prime }}x_{i_{j}}^{P_{j}(n)},$ where $i_{j}\in \{1,...,t\}$ for $
j=1,...,t^{\prime }$. If $G$ is a Lie group then by Ado's theorem one can
embed $G$ into a group of upper triangular matrices, in which setting the
matrix entries of both sides of the preceding formula are polynomials. It
follows that this formula is in fact valid with $n$ replaced by an arbitrary
real number $\eta $. \hfill {\usebox{\proofbox}}\vspace{11pt}

\begin{lemma}[Lazard]
\label{stewart-lem} There is a sequence of rational numbers $
\alpha_1,\beta_1,\dots,\alpha_m,\beta_m$, depending only on $s$, such that,
for all $x,y \in G$, 
\begin{equation*}
\exp(\log x + \log y) = x^{\alpha_1}y^{\beta_1} \dots
x^{\alpha_m}y^{\beta_m}.
\end{equation*}
Similarly there is a sequence of rational numbers $\gamma_1,\delta_1,\dots,
\gamma_m,\delta_m$ such that, for all $x,y \in G$, 
\begin{equation*}
\exp([\log x,\log y]) = x^{\gamma_1}y^{\delta_1} \dots
x^{\gamma_m}y^{\delta_m}.
\end{equation*}
\end{lemma}

\noindent\textit{Proof. } This is not quite the result stated in \cite
{stewart}, for example, where the expressions obtained are of the form $
w_1(x,y)^{\eta_1}w_2(x,y)^{\eta_2}\dots w_m(x,y)^{\eta_m}$ with each $w_i$ a
word in $x$ and $y$ (and $x^{-1}$, $y^{-1}$) and the $\eta_i$ rational
numbers. However it follows immediately from that formulation and Lemma \ref
{rat-powers}.\hfill{\usebox{\proofbox}}\vspace{11pt}

It follows immediately from the last two lemmas that there is a similar
expression for $\exp(\log x_1 + \dots + \log x_n)$. In order to get the
right bound for the last part of Theorem \ref{mainthm0} we will require a
certain amount of control over what this expression is.

\begin{lemma}[Expansion of sums]
\label{exp-sums} Let $G$ be an $s$-step nilpotent Lie group and suppose that 
$x_1,\dots,x_n \in G$, where $n \geq 2$. Then there is a sequence of
rational numbers $\alpha_1,\dots,\alpha_m$ depending on $n$, all of which
may be put over a common denominator of size bounded by $n^{C_s}$, and a
collection of indices $i_1,\dots,i_m \in \{1,\dots,n\}$ such that 
\begin{equation*}
\exp(\log x_1 + \dots + \log x_n) = x_{i_1}^{\alpha_1}\dots
x_{i_m}^{\alpha_m}.
\end{equation*}
\end{lemma}

\noindent\textit{Proof. } It clearly suffices to establish this when $n$ is
a power of two: if $n$ is not a power of two, let $n^{\prime r}$ be the
least power of two greater than $n$ and take $x_{n+1} = \dots = x_{n^{\prime
}} = \id_G$. To establish this case we proceed inductively, relating
the expansion for $n = k$ to that for $n = 2k$. We have 
\begin{equation*}
\exp(\log x_1 + \dots + \log x_k) = x_{i_1}^{\alpha_1} \dots
x_{i_m}^{\alpha_m}
\end{equation*}
and 
\begin{equation*}
\exp(\log x_{k+1} + \dots + \log x_{2k}) = x_{j_1}^{\alpha_1} \dots
x_{j_m}^{\alpha_m},
\end{equation*}
where the $\alpha_i$ may all be put over some denominator $Q$. Applying the
first of Lazard's expansions we may expand $\exp(\log x_1 + \dots + \log
x_{2k})$ as a product of terms of the preceding type, each to some rational
power over some fixed denominator $q_s$ depending only on $s$. Now expand
each of \emph{those} using Lemma \ref{rat-powers}: this results in an
expansion of $\exp(\log x_1 + \dots + \log x_{2k})$ as a product of terms $
x_{i_j}^{\beta_j}$, where all of the $\beta_j$ may be put over denominator $
Q q^{\prime }_s$ for some integer $q^{\prime }_s$ depending only on $s$. The
result follows immediately by induction. \hfill{\usebox{\proofbox}}\vspace{
11pt}

\section{Control by a nilbox}

In this section we prove Proposition \ref{key-prop-1}, the statement that an
approximate subgroup of an $s$-step nilpotent Lie group $G$ is controlled by
a nilbox. Recall that this, together with Proposition \ref{key-prop-2},
implies our main result.

Throughout this section $G$ is a simply-connected nilpotent Lie group and $
A\subseteq G$ is a $K$-approximate group. We write $\mathfrak{a}=\log
A=\{\log x:x\in A\}$ for the corresponding subset of the Lie algebra $
\mathfrak{g}$. Fisher, Katz and Peng \cite{fkp} used results close to those
of the previous section to prove that $\mathfrak{a}$ is close to invariant
under both addition and Lie bracket. We essentially recover this result in
Lemma \ref{a-abar} below. In their paper it was also remarked that one might
apply the Fre\u{\i}man-Ruzsa theorem in this setting, and we shall see how
this suggestion may be realised and used to prove Proposition \ref
{key-prop-1}. Some of the ideas here were also anticipated by the
foundational work of Tao \cite{tao-noncommutative}: in the last part of this
paper he described sets of small tripling in the Heisenberg group as being
precisely the sets that are roughly closed under both addition and
commutation, a fact he deduced by applying what amounts to the
Baker-Campbell-Hausdorff formula in this case.

\begin{lemma}[$\mathfrak{a}$ is almost invariant under nilcompletion]
\label{a-abar} Let $\overline{\mathfrak{a}}$ be the nilcompletion of $
\mathfrak{a}$. Then $|\mathfrak{a}+\overline{\mathfrak{a}}|\leq K^{C_{s}}| 
\mathfrak{a}|$. More generally $|m\overline{\mathfrak{a}}|\leq K^{mC_{s}}| 
\mathfrak{a}|$ for all $m\in \mathbb{N}$.
\end{lemma}

\noindent\textit{Proof. } It follows from Lemmas \ref{stewart-lem} and \ref
{rat-powers} that, for any $x_1,x_2,\dots \in G$ we have 
\begin{equation*}
\exp(\log x_1 + \log x_2 + [\log x_3,\log x_4] + [\log x_5,[\log x_6,\log
x_7]] + \dots) = x_1^{\eta_1} x_2^{\eta_2}\dots x_M^{\eta_M}
\end{equation*}
for some rationals $\eta_1,\eta_2,\dots,\eta_M$ and $M = C_s$. Here, all of
the commutators appearing in the definition of nilcompletion are featured.
Choose some integer $Q = C_s$ such that all of the rationals $\eta_i$ may be
put over the common denominator $Q$, and set $B :=
\{a^Q,a^{Q^2},\dots,a^{Q^s} : a\in A\}$ and $\mathfrak{b} = \log B$. Then
every element of $\exp(\mathfrak{b} + \overline{\mathfrak{b}})$ lies in $
A^{C_s}$, and hence by the iterated product set estimate (cf. Proposition 
\ref{approx-gp-basics}) we have $|\mathfrak{b} + \overline{\mathfrak{b}}|
\leq K^{C_s}|A|$.

Now suppose that $x\in \mathfrak{a}$. Then, noting that $\log (g^{t})=t\log
g $, we see that all of $Qx,Q^{2}x,\dots ,Q^{s}x$ lie in $\mathfrak{b}$.
Hence if $x_{1},x_{2},x_{3},\dots \in \mathfrak{a}$ then 
\begin{align*}
Q^{s}\cdot & \big(x_{1}+x_{2}+[x_{3}+x_{4}]+[x_{5},[x_{6},x_{7}]]+\dots 
\big) = \\
&
Q^{s}x_{1}+Q^{s}x_{2}+[Q^{s-1}x_{3},Qx_{4}]+[Q^{s-2}x_{5},[Qx_{6},Qx_{7}]]+
\dots
\end{align*}
lies in $\mathfrak{b}+\overline{\mathfrak{b}}$, that is to say $Q^{s}\cdot ( 
\mathfrak{a}+\overline{\mathfrak{a}})\subseteq (\mathfrak{b}+\overline{ 
\mathfrak{b}})$. The dilation map $Q^{s}:\mathfrak{g}\rightarrow \mathfrak{g}
$ is, of course, a bijection and so the result follows immediately. The last
claim follows immediately from the Ruzsa triangle inequality (\cite{tao-vu}
(2.6)) and its associated sum-set estimates (\cite{tao-vu} Cor. 2.23).\hfill 
{\usebox{\proofbox}}\vspace{11pt}

At this point we apply the Fre{\u{\i}}man-Ruzsa theorem, stated as Theorem 
\ref{frei-ruz} in the introduction. We will need a slightly stronger version
of this theorem than is commonly stated in the literature.

\begin{frei-ruz-again}
Suppose that $X$ is a $K$-approximate subgroup of $\mathbb{R}^{m}$. Then
there is a progression 
\begin{equation*}
P=\{l_{1}x_{1}+\dots +l_{k}x_{k}:|l_{i}|\leq L_{i}\}
\end{equation*}
with $k\leq K^{C}$ such that $X\subseteq P\subseteq K^{C}X$.
\end{frei-ruz-again}

\emph{Remarks on the proof.} There are two slight novelties here. The first
is minor and it is that the Fre{\u{\i}}man-Ruzsa theorem is normally only
stated for subsets of $\mathbb{Z}$, not of $\mathbb{R}^{m}$. This is
certainly addressed by the more general result of the second author and
Ruzsa \cite{green-ruzsa}, which is valid in an arbitrary abelian group (very
likely a more direct reduction to the $\mathbb{Z}$-case is also possible).
More seriously, the containment $P\subseteq K^{C}X$ is not normally stated
as part of the theorem and indeed we know of no reference in the literature
where it is explicitly mentioned. It may, however, be read without
difficulty out of \cite{green-ruzsa}. There one finds a proof that $2X-2X$
contains a progression $P_{0}$ of size at least $\exp (-K^{C})|X|$ and
dimension at most $K^{C}$. This is a standard ingredient in ``Ruzsa-style''
proofs of the Fre{\u{\i}}man-Ruzsa theorem and may also be found in Ruzsa's
original paper \cite{ruzsa-freiman} for subsets of $\mathbb{Z}$. To proceed
from such a statement to the Fre{\u{\i}}man-Ruzsa theorem one applies a
\textquotedblleft covering lemma\textquotedblright , and the most efficient
one in this context is implicit in Chang \cite{chang-freiman}. It is
explicitly stated in \cite[ Lemma 5.31]{tao-vu}, and using that result one
obtains 
\begin{equation*}
X\subseteq P_{0}-P_{0}+\{-1,0,1\}^{d}\cdot (v_{1},\dots ,v_{d})+x_{0},
\end{equation*}
for some $x_{0}\in X$, $v_{1},\dots ,v_{d}\in X-X=2X$ and where $d\leq K^{C}$
. Let $P$ be the set on the right. It is a progression of dimension at most $
K^{C}$, and it is fairly clearly contained in $(9+2d)X$. Also observe that
the $x_{i} $'s all lie in $4X.$\hfill {\usebox{\proofbox}}\vspace{11pt}

Putting the tools we have assembled so far together, we obtain the next
result, which is the key to our main result.

\begin{corollary}
\label{cor4.3} Suppose that $A \subseteq G$ is a $K$-approximate group and
let $\mathfrak{a} = \log A$. Then $\mathfrak{a}$ is contained in a \textup{
progression} $\mathfrak{p}$ of dimension at most $K^{C_s}$ whose
nilcompletion $\overline{\mathfrak{p}}$ is contained in $K^{C_s}\overline{
\mathfrak{\mathfrak{a}}}$.
\end{corollary}

\noindent \textit{Proof. } Recall Lemma \ref{a-abar}. Since $\mathfrak{a}
\subseteq \overline{\mathfrak{a}}$, the upper bound on $|\mathfrak{a}+ 
\overline{\mathfrak{a}}|$ certainly tells us that $|\mathfrak{a}+\mathfrak{a}
|\leq K^{C_{s}}|\mathfrak{a}|$. Applying Chang's version of the Fre{\u{\i}}
man-Ruzsa theorem we obtain a progression $\mathfrak{p}\subseteq K^{C_{s}}
\mathfrak{a }$ of dimension at most $K^{C_{s}}$ such that $\mathfrak{a}
\subseteq \mathfrak{p}$. The fact that $\overline{ \mathfrak{p}}\subseteq
K^{C_{s}}\overline{\mathfrak{a}}$ follows from Lemma \ref{lem5.1} (i)
.\hfill {\usebox{\proofbox}}\vspace{
11pt}

Now given the progression $\mathfrak{p} = \{l_1 x_1 + \dots + l_k x_k :
|l_i| \leq L_i\}$ it is very tempting to consider the nilbox $\mathfrak{B}
(x_1,\dots,x_k;L)$. One cannot quite use this to prove Proposition \ref
{key-prop-1}, however, since it is not necessarily the case that $\exp(x_i)
\in \langle A \rangle$ for $i = k+1,\dots,t$. This may be rectified by the
simple expedient of taking suitable dilates $M \cdot x_i$ of the generators $
x_i$, $i=1,\dots,k$. For $i=k+1,\dots,t$ write $x_i^{(M)} := M^{\chi(i)}
\cdot x_i$ and note that there is a Lie algebra homomorphism $\pi_M : 
\mathfrak{n}_{k,s} \rightarrow \mathfrak{g}$ such that $\pi_M(x_i) =
x_i^{(M)}$ for $i=1,\dots,t$ (this being the one induced by mapping $X_i$ to 
$M \cdot x_i$ for $i=1,\dots,k$).

\begin{lemma}
\label{lem4.6} Suppose that $A \subseteq G$ and that $\mathfrak{a} := \log A$
. Suppose that $\mathfrak{p}$ is a generalised progression as above, that $
\mathfrak{p}$ contains $\mathfrak{a}$ and that $\overline{\mathfrak{p}}$ is
contained in $m\overline{\mathfrak{a}}$. Let $M \geq 1$ be an integer. Then

\begin{enumerate}
\item The nilbox $\mathfrak{B}(x_1,\dots,x_k;L)$ contains $\mathfrak{a}$ and
is contained in $(4k)^s m \overline{\mathfrak{a}}$;

\item There is some $M \leq (2mk)^{C_s}$ for which the elements $
\exp(x_i^{(M)})$, $i=1,\dots,t$, all lie in the group $\langle A\rangle$
generated by $A$.
\end{enumerate}
\end{lemma}

\noindent\textit{Proof. } The first part of (i) is totally obvious since $
\mathfrak{p}$ contains $\mathfrak{a}$. To establish the second part we
employ an argument identical to the one in the proof of Proposition \ref
{prop33} to obtain the inclusions $\mathfrak{B}(x_1,\dots,x_k;L) \subseteq
(4k)^{s}\overline{\mathfrak{p}} \subseteq (4k)^sm\overline{\mathfrak{a}}$,
as required.

Let us turn to (ii). We take $M$ to have the special form $R^s$ for some
integer $R \geq 1$ to be specified. Then for each $i = 1,\dots,t$ the
element $x_i^{(M)} = M^{\chi(i)} \cdot x_i$ may be written as the sum of at
most $4^s m$ elements of the form $M^r[a_1,[a_2,[\dots]]]$, where the $a_i$
lie in $\mathfrak{a}$ and $r \geq 1$. We may each such element as $
[b_1,[b_2,[\dots]]]$ where the commutator has the same shape and each $b_i$
is $R^{u_i}a_i$ for some $u_i \geq 1$.

Now we simply expand $\exp(x_i^{(M)})$ using Lemmas \ref{rat-powers}, \ref
{stewart-lem} and \ref{exp-sums}, obtaining a finite product $
y_{i_1}^{\eta_1}y_{i_2}^{\eta_2}\dots$ with $y_i := \exp(b_i)$ in which, it
may be confirmed, the rationals $\eta_j$ may all be put over some common
denominator of size at most $(2ks)^{C_s}$. Taking $R$ to be this common
denominator it follows that each $y_{i_j}^{\eta_j}$ lies in the group
generated by $A$, which is what we wanted to prove.\hfill{\usebox{\proofbox}}
\vspace{11pt}

We are at last in a position to conclude the proof of Proposition \ref
{key-prop-1}, at least given a small result on nilboxes and their dilations
from the next section.\vspace{11pt}

\emph{Proof of Proposition \ref{key-prop-1}.} As ever $A \subseteq G$ is a $K
$-approximate group and $\mathfrak{a} := \log A$. Applying Corollary \ref
{cor4.3} we obtain a progression $\mathfrak{p}$ of dimension $K^{C_s}$ whose
nilcompletion $\overline{\mathfrak{p}}$ is contained in $K^{C_s}\overline{
\mathfrak{a}}$. Applying Lemma \ref{lem4.6} (i) we obtain a nilbox $
\mathfrak{B}(x_1,\dots,x_k;L)$ of dimension $k \leq K^{C_s}$ which contains $
\mathfrak{a}$ and is contained in $K^{C_s}\overline{\mathfrak{a}}$. It
follows from this last inclusion and Lemma \ref{a-abar} that $|\mathfrak{B}
(x_1,\dots,x_k;L)| \leq e^{K^{C_s}}|A|$, and hence $\exp(\mathfrak{B}
(x_1,\dots,x_k;L))$ $e^{K^{C_s}}$-controls $A$. Finally we may apply Lemma 
\ref{lem4.6} (ii) to find an $M = K^{C_s}$ such that the generators $
\exp(x_i^{(M)})$, $i = 1,\dots,t$, all lie in $\langle A\rangle$. By
Corollary \ref{dilation-compare}, $\exp(\mathfrak{B}(x_1^{(M)},
\dots,x_k^{(M)};L))$ also $e^{K^{C_s}}$-controls $A$. This concludes the
proof of Proposition \ref{key-prop-1}.\hfill{\usebox{\proofbox}}\vspace{11pt}

\section{Nilboxes, nilpotent progressions and control}

\label{nilbox}

To conclude the proof of our main theorem we must establish Proposition \ref
{key-prop-2}, which asserts that the exponential of a nilbox is efficiently
controlled by a nilpotent progression. We use the need for this as an excuse
to develop the relationship of nilboxes and nilpotent progressions more
generally.

Let us recall the statement of Proposition \ref{key-prop-2}.

\begin{key-prop-2-repeat}[Nilpotent progressions control nilboxes]
Suppose that $G$ is an $s$-step simply-connect nilpotent Lie group with Lie
algebra $\mathfrak{g}$. Suppose that $x_1,\dots,x_k \in G$, and write $u_i
:= \exp(x_i)$. Let $L = (L_1,\dots,L_k)$ be a vector of positive integer
lengths. Then the nilpotent progression $P(u_1,\dots,u_k;L)$ $e^{k^{C_s}}$
-controls $\exp(\mathfrak{B}(x_1,\dots,x_k;L))$.
\end{key-prop-2-repeat}

In fact we shall show that $\exp(\mathfrak{B}(x_1,\dots,x_k;L))$ $e^{k^{C_s}}
$-controls $P(u_1,\dots,u_k;L)$ as well. This notion of \emph{mutual control}
, where two sets $A$ and $B$ in some group $K$-control one another, is very
useful since if $\pi$ is a group homomorphism then $\pi(A)$ and $\pi(B)$
also $K$-control one another. Indeed the inclusions $A \subseteq (X \cdot B)
\cap (B \cdot X)$ and $B \subseteq (X^{\prime }\cdot A) \cap (A \cdot
X^{\prime })$ imply $\pi(A) \subseteq (\pi(X) \cdot \pi(B)) \cap (\pi(B)
\cdot \pi(X))$ and $\pi(B) \subseteq (\pi(X^{\prime }) \cdot \pi(A)) \cap
(\pi(A) \cdot \pi(X^{\prime }))$, conditions which automatically imply that $
|\pi(A)| \leq K|\pi(B)|$ and $|\pi(B)| \leq K|\pi(A)|$. In our setting, the
upshot of this together with the commutativity of the diagram 
\begin{equation}
\begin{CD} \mathfrak{n}_{k,s} @>\exp>> N_{k,s} \\ @VVV @VVV \\ \mathfrak{g}
@>\exp>> G\end{CD}
\end{equation}
is that we need only establish Proposition \ref{key-prop-2} in the free
setting.

Suppose then that $X_1,\dots,X_k,X_{k+1},\dots, X_t$ is an adapted basis for
the free $s$-step nilpotent Lie algebra $\mathfrak{n}_{k,s}$, as described
in Section \ref{strategy-sec}. Let us introduce the shorthand
\begin{equation*}
\mathfrak{B}(k,s;L) := \mathfrak{B}(X_1,\dots,X_k;L) := \{l_1 X_1 + \dots +
l_t X_t : |l_i| \leq L^{\chi(i)}\}
\end{equation*}
for the free nilbox on $k$ generators. We will require the variants 
\begin{equation*}
\mathfrak{B}(k,s;L,Q) := \{l_1 X_1 + \dots + l_t X_t : |l_i| \leq
L^{\chi(i)}, Q |l_i\}
\end{equation*}
for integers $Q \geq 1$. Write $u_i = \exp(X_i)$, $i = 1,\dots,k$ and
suppose that $u_{k+1},\dots,u_t$ is the ordered list of basic commutators in
the $u_i$ described in the introduction. We also introduce shorthands for
the nilpotent progressions based on the $u_i$, thus 
\begin{equation*}
P(k,s;L) := P(u_1,\dots,u_k;L) := \{ u_1^{l_1} \dots u_t^{l_t} : |l_i| \leq
L^{\chi(i)}\}.
\end{equation*}
Finally we introduce the variant 
\begin{equation*}
P(k,s;L,Q) := \{ u_1^{l_1} \dots u_t^{l_t} : |l_i| \leq L^{\chi(i)}, Q |
l_i\}.
\end{equation*}
The key proposition linking these two types of object is the following
proposition, which we will establish in the appendix after developing some
basic theory of ``coordinates'' on the free nilpotent group. Those results
ultimately rest on the Baker-Campbell-Hausdorff formula.

\begin{proposition}
\label{inclusions} There is some integer $Q_{s}$ and constants $c=c_{s}$, $
C=C_{s}$ such that for any positive integer $Q$ divisible by $Q_{s}$ we have 
\begin{equation*}
\exp(\mathfrak{B}(k,s;cL,Q))\subseteq P(k,s;L);
\end{equation*}
\begin{equation*}
P(k,s;cL,Q)\subseteq \exp(\mathfrak{B}(k,s;L))
\end{equation*}
and such that uniformly in $\rho ,\rho ^{\prime }<10$, we have 
\begin{equation*}
\exp(\mathfrak{B}(k,s;\rho L,Q))\cdot \exp(\mathfrak{B}(k,s;\rho^{\prime
}L,Q))\subseteq \exp(\mathfrak{B}(k,s; (\rho +\rho ^{\prime}+C\rho
\rho^{\prime})L;Q)
\end{equation*}
and 
\begin{equation*}
P(k,s;\rho L,Q)\cdot P(k,s;\rho^{\prime }L,Q)\subseteq P(k,s;(\rho
+\rho^{\prime }+C\rho \rho^{\prime})L,Q)).
\end{equation*}
In fact 
\begin{equation*}
P(k,s,\rho L) \cdot P(k,s,\rho^{\prime }L) \subseteq P(k,s;(\rho +
\rho^{\prime }+ C\rho \rho^{\prime })L).
\end{equation*}
\end{proposition}

\emph{Remark.} The final inclusion here easily implies that if we set $
X_{\rho} := P(k,s;(\rho - C^{\prime}\rho^2) L)$ for an appropriate $C^{\prime}$
then $X_{\rho}\cdot X_{\rho^{\prime }} \subseteq X_{\rho + \rho^{\prime }}$.
Furthermore Proposition \ref{control-balls} below together with a little
calculation implies that $X_{2\rho}$ is $\exp(C_sk^s)$-controlled by $
X_{\rho}$. These two facts imply that the system $(X_{\rho})_{\rho \leq 4}$
forms a \emph{Bourgain system} in the sense of \cite{green-sanders,sanders},
thereby providing a link between our work and that of Sanders.\vspace{11pt}

Proposition \ref{key-prop-2} follows immediately from Proposition \ref
{inclusions} and the next result.

\begin{proposition}
\label{control-balls} Let $k,s \geq 1$ be integers and let $Q,Q^{\prime
}\geq 1$ be two integers with $Q| Q^{\prime }$ and $Q_s |Q^{\prime }$, where 
$Q_s$ is the quantity appearing the previous proposition. Let $\lambda< 1$
be a positive real number. Then $\exp(\mathfrak{B}(k,s,L;Q))$ and $\exp(
\mathfrak{B}(k,s,\lambda L;Q^{\prime }))$ mutually $\gamma$-control one
another, as do $P(k,s,L;Q)$ and $P(k,s,\lambda L;Q^{\prime })$, where we may
take $\gamma = (1 + \frac{Q^{\prime }}{\lambda Q} )^{C_sk^s}$.
\end{proposition}

\noindent\textit{Proof. } In view of the third inclusion of Proposition \ref
{inclusions}, there is some $\lambda^{\prime }\geq c_s \lambda$ such that $
\exp(\mathfrak{B}(k,s,\lambda^{\prime }L;Q^{\prime 3}\subseteq \exp(
\mathfrak{B}(k,s,\lambda L,Q^{\prime }))$. Hence to conclude the proof it
suffices to show, in view of Proposition \ref{approx-gp-basics}, the bound 
\begin{equation*}
\frac{|\mathfrak{B}(k,s,\lambda^{\prime }L;Q^{\prime })|}{|\mathfrak{B}
(k,s,L;Q)|} \geq (1 + \frac{Q^{\prime }}{\lambda Q})^{-C_sk^s}.
\end{equation*}

However one has 
\begin{equation*}
|\mathfrak{B}(k,s,\lambda^{\prime }L;Q^{\prime })| = \prod_{j=1}^t (2\lfloor 
\frac{(\lambda^{\prime} L)^{\chi(j)}}{Q^{\prime }}\rfloor + 1) \geq \prod_{j=1}^t
(2 \lfloor \frac{\lambda^{\prime s}}{Q^{\prime }} L^{\chi(j)} \rfloor + 1)
\end{equation*}
whereas 
\begin{equation*}
|\mathfrak{B}(k,s,L;Q)| = \prod_{j=1}^t (2\lfloor \frac{L^{\chi(j)}}{Q}
\rfloor + 1) .
\end{equation*}
Note, however, the inequality 
\begin{equation*}
2\lfloor x \rfloor + 1 \leq 2x + 1 \leq \frac{2}{\alpha}(\lfloor \alpha
x\rfloor + 1) + 1 \leq (1 + \frac{1}{\alpha})(2\lfloor \alpha x\rfloor + 1),
\end{equation*}
valid for any $\alpha,x \geq 0$. It follows that 
\begin{equation*}
\frac{|\mathfrak{B}(k,s,\lambda^{\prime }L;Q^{\prime })|}{|\mathfrak{B}
(k,s,L;Q)|} \geq (1 + \frac{Q^{\prime }}{\lambda^{\prime }Q})^{-t}.
\end{equation*}
Finally note, as remarked in the introduction, that $t \leq (4k)^s$.

The proof for the nilprogressions $P(k,s;L)$ is essentially identical.\hfill{
\usebox{\proofbox}}\vspace{11pt}

Let us record a particular application of this that we used in \S \ref{sec4}.

\begin{corollary}[Control by dilated balls]
\label{dilation-compare} Let $G$ be a simply-connected $s$-step nilpotent
Lie group with Lie algebra $\mathfrak{g}$. Let $x_1,\dots,x_k \in \mathfrak{g
}$, and suppose that $M,M^{\prime }\geq 1$ are integers. Then $\exp(
\mathfrak{B}(Mx_1,\dots,Mx_k;L))$ and $\exp(\mathfrak{B}(M^{\prime
}x_1,\dots,M^{\prime }x_k;L))$ mutually $(MM')^{C_sk^s}$-control one
another.
\end{corollary}

\noindent\textit{Proof. } The introduction of $M^{\prime}$ is merely to
make the statement look symmetric. We only needed the corollary in the case $
M^{\prime }=1$, and in fact the general case clearly follows from this
special one using the transitivity of control. Suppose, then, that $
M^{\prime }= 1$. Since we are dealing with \emph{mutual} control, a notion
which persists under homomorphisms by the remarks at the beginning of the
section, it suffices to work in the free setting.

One may easily check the inclusions 
\begin{equation*}
\mathfrak{B}(MX_1,\dots,MX_k;M^{-s}L,Q)\subseteq \mathfrak{B}
(X_1,\dots,X_k;L,Q)
\end{equation*}
and 
\begin{equation*}
\mathfrak{B}(X_1,\dots, X_k;L,M^sQ)\subseteq \mathfrak{B}(MX_1,
\dots,MX_k;L,Q)
\end{equation*}
for any integer $Q \geq 1$. Thus 
\begin{equation*}
\mathfrak{B}(MX_1,\dots,MX_k;L) \supseteq \mathfrak{B}(X_1,\dots,X_k;L,M^s)
\supseteq \mathfrak{B}(MX_1,\dots,MX_k;M^{-s}L,M^s),
\end{equation*}
which implies by Proposition \ref{control-balls} that 
\begin{equation*}
\exp(\mathfrak{B}(MX_1,\dots,MX_k;L))
\end{equation*}
$M^{C_sk^s}$-controls 
\begin{equation*}
\exp(\mathfrak{B}(X_1,\dots,X_k;L,M^s)).
\end{equation*}
By another application of Proposition \ref{control-balls} this, in turn, $
M^{C_sk^s}$-controls $\exp(\mathfrak{B}(X_1,\dots,X_k;L))$. The inverse
relationship may be obtained very similarly. \hfill{\usebox{\proofbox}}
\vspace{11pt}

\appendix

\section{Fre\u{\i}man invariance of nilpotent progressions}

In this section we show that nilpotent progressions are preserved under 
\emph{Fre\u{\i}man homomorphisms} if the sidelengths $L$ are sufficiently
large. Thus, in a sense, they are well-defined multiplicative-combinatorial
objects. We recall the definition of Fre\u{\i}man homomorphism. Suppose that 
$A$ and $B$ are two sets in ambient groups and that $\phi : A \rightarrow B$
is a map. We say that $\phi$ is a Fre\u{\i}man $k$-homomorphism provided
that for all $a_1,\dots,a_k \in A$ and all choices of $\varepsilon_1,\dots,
\varepsilon_k \in \{-1,0,1\}$ the conditions $a_1^{\varepsilon_1} \dots
a_k^{\varepsilon_k} = \id$ implies that $\phi(a_1)^{\varepsilon_1}
\dots \phi(a_k)^{\varepsilon_k} = \id$.

\begin{proposition}[Invariance under Fre\u{\i}man isomorphism]
Suppose $G$ and $H$ are two groups and $s\in \mathbb{N}$. There is a
constant $C_{s}\geq 1$ such that the following holds. Let $A$ be an $s$-step
nilpotent progression with $k$ generators in $G$ and side lengths $
L=(L_{1},...,L_{k})$ with $L_{i}\geq C_{k,s}$ for $i=1,...,k.$ Let $\phi $
be a Fre\u{\i}man $3$-homomorphism from $A$ onto a subset $B$ of $H$. Then $B
$ is also an $s$-step nilpotent progression with $k$ generators.
\end{proposition}

\noindent\textit{Proof. } Every nilpotent progression has the form 
\begin{equation}  \label{nil-box-4}
A = \{u_1^{l_1} \dots u_t^{l_t} : |l_j| \leq L^{\chi(j)}\}
\end{equation}
for some elements $u_1,\dots,u_t$. However not every set of this form is a
nilpotent progression: indeed $u_{k+1},\dots,u_t$ are specific commutators
involving $u_1,\dots,u_k$, and in particular there are words $w_i$, $
i=k+1,\dots,t$, independent of the underlying group, such that $u_i =
w_i(u_1,\dots,u_k) $. Moreover it is not hard to see that \emph{if} the $u_i$
satisfy these conditions and if all $(s+1)$-fold commutators of the $u_i$
equal the identity then the object \eqref{nil-box-4} does define a nilpotent
progression.

Now if the lengths $L = (L_1,\dots,L_k)$ are sufficiently great then all
initial segments of all these words $w_j$ lie in the nilpotent progression $
P(u_1,\dots,u_k;L)$. This follows from the fifth inclusion of Proposition 
\ref{inclusions}, or else it may be verified more explicitly by taking each
initial segment $u_{i_1},\dots,u_{i_m}$ and commutating until all copies of $
u_1$ are at the left, then repeating this process for $u_2$ and so on.

Now let $\phi : A \rightarrow B$ be a map. Observe that if $a_{1},...,a_{m}$
are elements of $A$ or $A^{-1}$ such that each initial segment $a_{1}\cdots
a_{i}$ belongs to $A$ for all $i=1,\dots,m$ then $\phi (a_{1}\cdots
a_{m})=\phi (a_{1})\cdots\phi (a_{m})$. Also, by an easy induction, we have $
\phi(a^j) = \phi(a)^j$ whenever all the powers $a,a^2,\dots,a^j$ lie in $A$.
Writing $v_i := \phi(u_i)$ for $i=1,\dots,t$, it follows from these
observations and the analysis of the preceding paragraph that the $v_i$
satisfy the same words $v_i = w_{i}(v_1,\dots,v_k)$, and also that all $(s+1)
$-fold commutators of the $v_i$ equal the identity.

It follows that 
\begin{equation*}
\{ v_1^{l_1} \dots v_t^{l_t} : |l_j| \leq L^{\chi(j)}\}
\end{equation*}
is a nilpotent progression in $H$. Furthermore several more applications of
the observation we made in the last paragraph confirm that 
\begin{equation*}
\phi(u_1^{l_1} \dots u_t^{l_t}) = v_1^{l_1} \dots v_t^{l_t},
\end{equation*}
and so this nilpotent progression is precisely equal to $B$.\hfill{
\usebox{\proofbox}}\vspace{11pt}

\section{On coordinates in the free nilpotent Lie group}

\label{coord-append}

Throughout this appendix we will be working in the free $s$-step nilpotent
Lie algebra $\mathfrak{n}_{k,s}$ on $k$ generators and with the
corresponding free nilpotent Lie group $N_{k,s}$. We suppose that an adapted
basis (see the introduction for definitions) $X_1,\dots,X_t$ for $\mathfrak{n
}_{k,s}$ has been chosen, and that $u_1,\dots,u_t$ is the corresponding
ordered list of group commutators in $N_{k,s}$. We will also use the weight
function $\chi :[t] \rightarrow \mathbb{N }_0^k$ introduced in the
introduction.

Our aim in this section is to establish Proposition \ref{inclusions}, which
the reader may care to recall now. We shall be quite brief in our treatment,
which depends on a study of \emph{coordinates} in the following sense.

\begin{definition}[Coordinates]
Suppose that $x \in N_{k,s}$. Then we define the group coordinates $\psi_{
\gp}(x)$ to be $(x_1,\dots,x_t)$, where $x_1,\dots,x_t$ are the unique
complex numbers such that $x = u_1^{x_1} \dots u_t^{x_t}$. We define the
algebra coordinates $\psi_{\alg}(x)$ to be $(x^{\prime
}_1,\dots,x^{\prime }_t)$, where $x^{\prime }_1,\dots,x^{\prime }_t$ are the
unique complex numbers such that $x = \exp(x^{\prime }_1 X_1 + \dots +
x^{\prime }_t X_t)$.
\end{definition}

\emph{Remarks.} The existence and uniqueness of the group coordinates is not
obvious, and it will be a by-product of our analysis. The algebra
coordinates are the same thing as \textquotedblleft exponential coordinates
of type I\textquotedblright\ as featured for example in \cite
{corwin-greenleaf} (hence the notation, with a single dash). The group
coordinates are \emph{not} quite the same thing a exponential coordinates of
type II. If $x\in N_{k,s}$ then to find the type II coordinates $
(x_{1}^{\prime \prime },\dots ,x_{t}^{\prime \prime })$ one expresses $x$ as 
$\exp (x_{1}^{\prime \prime }X_{1})\dots \exp (x_{t}^{\prime \prime }X_{t})$
. We will encounter type II coordinates again in a short while. \vspace{11pt}

Everything will follow from the Baker-Campbell-Hausdorff formula, which
states that 
\begin{equation*}
\exp(X)\exp(Y) = \exp(X + Y + \frac{1}{2}[X,Y] + \frac{1}{12}[X,[X,Y]] +
\dots)
\end{equation*}
It is not important to know what the rationals numbers here are, and indeed
they are rather complicated to describe. All that is important is that the
series on the right is finite in an $s$-step nilpotent group, and that all
the rationals occurring have complexity $O_s(1)$.

With a little thought, this leads directly to the following description of
multiplication in algebra coordinates.

\begin{lemma}[Multiplication in algebra coordinates]
\label{typei-coordinates} Suppose that $x,y \in N_{k,s}$ and that $\psi_{ 
\alg}(x) = (x_1,\dots,x_t)$ whilst $\psi_{\alg}(y) =
(y_1,\dots,y_t)$ . Then 
\begin{equation*}
\psi_{\alg}(xy) = (P_1(x_i,y_i),\dots,P_t(x_i,y_i)).
\end{equation*}
Here, each $P_{j}$ is a polynomial of the form $x_{j} + y_j +
\sum_{\alpha,\beta \neq 0} C^{(j)}_{\alpha,\beta } x^{\alpha}y^{\beta}$,
where the $C^{(j)}_{\alpha,\beta}$ are rationals with complexity $O_s(1)$, $
\alpha = (\alpha_1,\dots,\alpha_t)$ and $\beta = (\beta_1,\dots,\beta_t)\in 
\mathbb{N} _0^t$ are multi-indices, and $x^{\alpha}$ means $
x_1^{\alpha_1}\dots x_t^{\alpha_t}$. Furthermore $C^{(j)}_{\alpha,\beta}$ is
only nonzero if $\sum_{l \in [t]} \chi(l) (\alpha_l + \beta_l) = \chi(j)$.
\end{lemma}

The third inclusion of Proposition \ref{inclusions} follows quickly from
this, taking $Q_s$ to be the least common multiple of the denominators of
all the $C^{(j)}_{\alpha,\beta}$.

We may also say something about the transformation which takes the group
coordinates of a point $x \in N_{k,s}$ and outputs the algebra coordinates.
Before doing this it is convenient to set up a notion of \emph{degree}.
Suppose that $(x_1,\dots,x_t) $ are variables, to be thought of as
coordinates. Now let $z_1,\dots,z_k$ be further variables, and substitute $
x_j = z^{\chi(j)} = z_1^{\chi(j)_1} \dots z_k^{\chi(j)_k}$. When we speak of
the \emph{degree} $\deg(P)$ of a polynomial $P = P(x_1,\dots,x_t)$ involving
the $x_i$, we shall mean the total degree in the $z_i$ of $P(z_1,\dots,z_k)$
after this substitution has been made. We shall also write $\deg_i$ for the
degree with respect to $z_i$.

\begin{definition}[Privileged coordinate change]
Let $\phi : \mathbb{C}^t \rightarrow \mathbb{C}^t$ be a polynomial map. We
say that $\phi$ is a \emph{privileged coordinate change} if $(\phi(x))_j =
x_j + P_j(x_1,\dots,x_t)$, where $\deg_i(P_j(x)) \leq \deg_i(x_j)$ for all $
i = 1,\dots,k$ and $P_j$ depends only on those variables $x_l$ with $
\deg(x_l) < \deg(x_j)$.
\end{definition}

The inverse of a privileged coordinate change is another privileged
coordinate change, as is the composition of two such coordinate changes. We
leave the proof as an exercise. If all the coefficients of the $P_j$ are
rationals with complexity $O_s(1)$ then we say that $\phi$ is of \emph{
bounded complexity}; the inverse and composition of privileged coordinate
changes of bounded complexity are also privileged coordinate changes of
bounded complexity, albeit with worsenings of the unspecified constants $
O_s(1)$.

\begin{lemma}[Group coordinates to algebra coordinates]
\label{lemma-i-ii} Suppose that $x \in N_{k,s}$ and that $\psi_{\gp
}(x) = (x_1,\dots,x_t)$ and $\psi_{\alg}(x) = (x^{\prime
}_1,\dots,x^{\prime }_t)$. Then the mapping $(x_j) \rightarrow (x^{\prime
}_j)$ and its inverse are privileged coordinate changes of bounded
complexity.
\end{lemma}

\noindent\textit{Proof. } Write $Y_i := \log (u_i)$ for $i = 1,\dots,t$. We
claim that the relation between the $Y_i$ and the $X_i$ is a rather special
one: we have 
\begin{equation*}
Y_j = X_j + \sum_{m} \mu_{jm} X_m
\end{equation*}
where the sum is restricted to those $m$ for which $\chi(m) \geq \chi(j)$
pointwise but $\chi(m) \neq \chi(j)$. This may be established by induction
on the order of the commutator $u_i$ using the Baker-Campbell-Hausdorff
formula: we leave the details to the reader. Suppose that $\psi_{\alg
}(x) = (x^{\prime }_1,\dots,x^{\prime }_t)$, that is to say $\log x =
x^{\prime }_1 X_1 + \dots + x^{\prime }_t X_t$. We may also represent $\log x
$ as $y_1 Y_1 + \dots + y_t Y_t$ in ``exponential coordinates of type I''
relative to the basis $Y_i$, and it is not hard to check that the change of
coordinates map $(x^{\prime }_i) \rightarrow (y_i)$ is a privileged
coordinate change of bounded complexity.

In view of the group closure properties of the notion of privileged
coordinate change it suffices to show that the map $(x_i) \rightarrow (y_i)$
is a privileged coordinate change of bounded complexity. Note that $x =
\exp(x_1 Y_1) \dots \exp(x_t Y_t)$ (so the group coordinates $\psi_{\gp
}(x) = (x_i)$ are actually the same thing as the ``exponential coordinates
of type II'' relative to the basis $Y_1,\dots,Y_t$). The desired property
follows from repeated application of the Baker-Campbell-Hausdorff formula:
once again we leave the precise details to the reader.\hfill{
\usebox{\proofbox}}\vspace{11pt}

\emph{Remark.} The existence and uniqueness of group coordinates follows
from (the proof of) this lemma.\vspace{11pt}

The first and second inclusions of Proposition \ref{inclusions} follow very
quickly from this lemma. To establish the fourth and fifth inclusions of
that proposition, it suffices to prove the following ``group'' variant of
Lemma \ref{typei-coordinates}.

\begin{lemma}[Multiplication in group coordinates]
\label{typeii-coordinates} Suppose that $x,y \in N_{k,s}$ and that $
\psi_{gp}(x) = (x_1,\dots,x_t)$ whilst $\psi_{\gp}(y) = (y_1,\dots,y_t)
$ . Then 
\begin{equation*}
\psi_{\gp}(xy) = (P_1(x_i,y_i),\dots,P_t(x_i,y_i)).
\end{equation*}
Here, each $P_{j}$ is a polynomial mapping $\mathbb{Z}^t \times \mathbb{Z}^t$
to $\mathbb{Z}$ of the form $x_{j} + y_j + \sum_{\alpha,\beta \neq 0}
C^{(j)}_{\alpha,\beta } x^{\alpha}y^{\beta}$, where the $C^{(j)}_{\alpha,
\beta}$ are rationals with complexity $O_s(1)$, $\alpha =
(\alpha_1,\dots,\alpha_t)$ and $\beta = (\beta_1,\dots,\beta_t)\in \mathbb{N}
_0^t$ are multi-indices, and $x^{\alpha}$ means $x_1^{\alpha_1}\dots
x_t^{\alpha_t}$. Furthermore $C^{(j)}_{\alpha,\beta}$ is only nonzero if $
\sum_{l \in [t]} \chi(l) (\alpha_l + \beta_l) = \chi(j)$.
\end{lemma}

\noindent\textit{Proof. } Combine Lemma \ref{lemma-i-ii} with Lemma \ref
{typei-coordinates} to conclude that Lemma \ref{typei-coordinates}. The fact
that each $P_i$ maps $\mathbb{Z}^t \times \mathbb{Z}^t$ to $\mathbb{Z}$
follows from the fact that the set $\{u_1^{n_1} \dots u_t^{n_t} :
n_1,\dots,n_t \in \mathbb{Z}\}$ is a group (the free $s$-step nilpotent
group $\Gamma_{k,s}$ on $k$ generators). This may be verified by repeated
commutation, taking the product of two such elements and moving all copies
of $u_1$ to the left, then all copies of $u_2$, and so on. \hfill{
\usebox{\proofbox}}\vspace{11pt}

\emph{Remark.} Similar issues to those addressed by the last lemma are
discussed in Leibman \cite{leibman}.

\end{document}